\documentclass[doublespacing]{elsart3-1}

\usepackage{amsfonts,amsmath,eucal,amssymb}
\usepackage[english,francais]{babel}

\usepackage[T1]{fontenc}
\usepackage[latin1]{inputenc}

\usepackage{url}
\usepackage{enumerate}

\usepackage[dvips]{graphicx}
\usepackage[all,2cell,ps]{xy}

\newcommand{\lulab}[1]{\ar@{}[l]_<<{#1}} 
\newcommand{\rulab}[1]{\ar@{}[r]^<<{#1}}
\newcommand{\ldlab}[1]{\ar@{}[l]^<<{#1}}
\newcommand{\rdlab}[1]{\ar@{}[r]_<<{#1}}
\newcommand{\edge}[1]{\ar@{-}[#1]} 
\newcommand{\node}{*+[o][F-]{ }} 

\newtheorem{theorem}{Theorem}[section]

\newtheorem{e-proposition}[theorem]{Proposition}

\newtheorem{e-definition}[theorem]{Definition} 

\newtheorem{theoreme}{Th\'eor\`eme}[section]

\newtheorem{definition}[theoreme]{D\'efinition\rm}

\renewcommand{\theequation}{\arabic{equation}}
\newcommand{\halmos}{\rule{1ex}{1.4ex}}
\newenvironment{proof}{\noindent {\em Proof}.\ \ %
}{\hspace*{\fill}$\halmos$\medskip}

\newcommand{\poset}{\ensuremath{{\mathcal{S}}}}

\def\og{\leavevmode\raise.3ex\hbox{$\scriptscriptstyle\langle\!\langle$~}}
\def\fg{\leavevmode\raise.3ex\hbox{~$\!\scriptscriptstyle\,\rangle\!\rangle$}}

\newcommand{\cal}[1]{\ensuremath{\mathcal #1}}

\newcommand{\be}[1]{\begin{equation}\label{#1}}
 \newcommand{\ee}{\end{equation}}

\newcommand{\bt}[1]{\begin{theorem}\label{#1}}
 \newcommand{\bd}[1]{\begin{e-definition}\label{#1}}
 \newcommand{\bp}[1]{\begin{e-proposition}\label{#1}}

 \newcommand{\bdf}[1]{\begin{definition}\label{#1}}
 \newcommand{\edf}{\end{definition}}

 \newcommand{\et}{\end{theorem}}
\newcommand{\ed}{\end{e-definition}}
 \newcommand{\ep}{\end{e-proposition}}

\newcommand{\bpr}{\begin{proof}}
\newcommand{\epr}{\end{proof}}

 \newcommand{\bi}{\begin{itemize}}
 \newcommand{\ei}{\end{itemize}}
 \newcommand{\ben}{\begin{enumerate}}
 \newcommand{\een}{\end{enumerate}}

\def \Z {{\mathbb Z}}
\def \R {{\mathbb R}}

 \def \ra {\rightarrow }

 \def\L{\Lambda}

 \def\G{\Gamma}

\newcommand{\impt}[1]{{\it #1}}

\date{2006}
\begin{document}

\begin{frontmatter}

\selectlanguage{english}
\title{Monotonicity and complete monotonicity \\ for continuous-time Markov chains}

\vspace{-2.6cm}

\selectlanguage{francais}
\title{Monotonie et monotonie complète \\ des cha{\^i}nes de Markov à temps continu}

\selectlanguage{english}
 \author[Padoue]{Paolo {Dai Pra}}, 
 \ead{daipra@math.unipd.it}
 \author[Potsdam]{Pierre-Yves {Louis}},
\ead{louis@math.uni-potsdam.de}
 \author[Padoue]{Ida {Minelli}}
 \ead{minelli@math.unipd.it}
\vskip 0.5\baselineskip

 \address[Padoue]{Dipartimento di Matematica Pura e Applicata, Universit\`a di Padova, %
 Via Belzoni 7, 35131 Padova, Italy}
 \address[Potsdam]{Institut für Mathematik, Potsdam Universit\"at,  Am neuen Palais,10 -- Sans Souci, %
  D-14 415 Potsdam, Germany}

  \url{http://www.sciencedirect.com/science/article/pii/S1631073X06001646}

  Post-print of Comptes Rendus Mathematique,\\ Volume 342, Issue 12, 15 June 2006, pp. 965-970, ISSN 1631-073X \\
doi:10.1016/j.crma.2006.04.007\\
Creative Commons Attribution Non-Commercial No Derivatives License

\begin{abstract}
We analyze the notions of monotonicity and complete monotonicity for Markov Chains in continuous-time,
taking values in a finite partially ordered set. Similarly to what happens in discrete-time, the
two notions are not equivalent. However, we show that there are partially
ordered sets for which monotonicity and complete monotonicity coincide in continuous time but not in discrete-time.
{\it To cite this article: P. Dai Pra et al., C. R. Acad. Sci. Paris, Ser. I 342 (2006).}
\vskip 0.5\baselineskip

\selectlanguage{francais}
\noindent{\bf R\'esum\'e}
\vskip 0.5\baselineskip
\noindent
Nous étudions les notions de monotonie et de monotonie complète pour les processus de Markov (ou chaînes de Markov
à temps continu)
prenant leurs valeurs dans
un espace partiellement ordonné. Ces deux notions ne sont pas équivalentes, comme c'est le cas lorsque le temps est discret.
Cependant, nous établissons que pour certains ensembles partiellement ordonnés, l'équivalence a lieu en temps continu
bien que n'étant pas vraie en temps discret.
{\it Pour citer cet article~:P. Dai Pra et al., C. R. Acad. Sci. Paris, Ser. I 342 (2006)}
\end{abstract}
\end{frontmatter}

\selectlanguage{francais}
\section*{Version fran\c{c}aise abr\'eg\'ee}
L'utilisation des chaînes de Markov dans le cadre des algorithmes MCMC soulève de nombreuses questions au sein desquelles
la monotonie joue un rôle important. Deux notions de monotonie sont considérées pour les chaînes de Markov à valeurs dans
un espace \impt{partiellement ordonné} $(S,<)$ (\impt{poset} selon la terminologie anglaise).
Nous supposons $S$ fini et les chaînes (resp. processus) de Markov homogènes en temps.
Dans la définition~\ref{mon1} nous reformulons la notion de \impt{monotonie}, équivalente à la définition usuelle via la propriété de
 stabilité des fonctions croissantes~$f$ sous l'action de l'\impt{opérateur de transition}
$T_t f(y) := \sum_{x \in S} f(x) P_t(x,y)$.

Pour différentes applications et en particulier à des fins de simulation,
une notion plus forte, dite de \impt{monotonie complète} a été introduite dans la Définition~\ref{mon2}.
En d'autres termes, cette {monotonie complète} signifie que l'on peut coupler simultanément
tous les processus partant de toutes les conditions initiales possibles, en préservant l'ordre entre ces processus.
Cette propriété s'avère essentielle pour l'utilisation de l'algorithme de Propp \& Wilson~(cf.~\cite{ProppWilson}).

Aucune méthode générale  n'existe pour caractériser simplement la monotonie complète à partir des
probabilités de transition ou du générateur infinitésimal.
De manière évidente, pour une chaîne (resp. processus )
de Markov donnée à valeurs dans~$S$, la monotonie complète implique la monotonie.
La question est alors de caractériser les posets~$S$ pour lesquels l'implication réciproque est vraie pour tout processus.
Pour les chaînes de Markov, {\it i.e.} lorsque le temps est discret,
ce problème a été résolu dans  \cite{FillMachida} (cf.~Theorem~\ref{fill-machida}
de la section~\ref{int}).
Dans cette note nous considérons cette question pour des processus de Markov (ou chaînes de Markov à temps continu) réguliers.

Soient $L = (L_{x,y})_{x,y \in S}$ le générateur infinitésimal d'un processus de Markov régulier et \\
\mbox{$S_2 := S \times S \setminus \{(x,x): x \in S\}$}.
Le générateur~$L$ peut être identifié à un élément du cône $(\R^+)^{S_2}$.

 Notre technique principale consiste à caractériser l'ensemble des générateurs monotones~${\cal{G}}_{mon}$
et celui des générateurs complètement
monotones~${\cal{G}}_{c.mon}$  comme des cônes de  $(\R^+)^{S_2}$.

La monotonie d'un générateur~$L$ est reformulée de
 manière géométrique~(\ref{Cone_mon_H-representation})
 dans la {\bf proposition~\ref{mon-cone}} en fonction des vecteurs~$W^{\G,x,y} \in \R^{S_2}$ (définis par~(\ref{vectori_W}))
où $x,y \in S$, $\G$ up-set (cf.~Définition~\ref{upset}).
 La représentation (\ref{Cone_mon_H-representation}) du cône~${\cal{G}}_{mon}$  comme intersection de demi-espaces
est appelée \impt{H-représentation} (cf.~\cite{FukudaProdon}).

Dans la {\bf proposition~\ref{prop2.3}} la monotonie complète d'un générateur~$L$ est caractérisée par sa
décomposition~(\ref{Cone_compl_mon_V-representation})
en tant que combinaison linéraire à coefficients positifs des vecteurs  $\mathbb{I}_f \in (\R^+)^{S_2}$ définis en~(\ref{vectori_I})
où $f \in {\cal{M}}$, {\it i.e.} où $f$ est une fonction croissante de $S$ dans $S$.
La représentation~(\ref{Cone_compl_mon_V-representation}) du cône~${\cal{G}}_{c.mon}$ via les rayons extrémaux $\mathbb{I}_f$
est appelée \impt{V-représentation}.

Puisque, pour toute fonction~$f$ croissante, $\G$ up-set, $x,y \in S$, on a
$ \langle \mathbb{I}_f,W^{\G,x,y} \rangle \geq 0$,
 on a bien l'inclusion ${\cal{G}}_{c.mon}\subseteq {\cal{G}}_{mon}$.
 Il apparaît dans un premier temps~({\bf Proposition~\ref{impl}})
que si l'équivalence entre les deux notions de monotonie a lieu pour toute chaîne de Markov à valeurs dans un
poset~$S$, alors cette équivalence demeure pour le poset~$S$ lorsque l'on considère des processus à temps continu.
La réciproque s'avère cependant fausse comme les résultats ci-après le prouve.

 Nous donnons la liste complète des posets
d'au plus cinq élements pour lesquels l'équivalence est fausse.
Nous n'avons malheureusement
pas réussi à caractériser tous les posets pour lesquels l'équivalence entre les deux notions de monotonie
est vérifiée, {\it i.e.} pour lesquels
${\cal{G}}_{c.mon}={\cal{G}}_{mon}$  a lieu.
Nous prouvons~:
\medskip

{\bf Proposition~\ref{Prop3.1}}
Les seuls posets~$S$ avec $\sharp S\leq 5$ tels que les deux notions de monotonie ne sont pas équivalentes
sont ceux dont le diagramme de Hasse (cf. définition {\it in} section~2~\cite{FillMachida}, p.~943)
est représenté dans la figure~\ref{5posets}.
\medskip

En particulier, cela signifie que ${\cal{G}}_{mon}={\cal{G}}_{c.mon}$ pour les deux posets importants
de cardinal quatre, à savoir le diamant et le noeud-papillon, dont les
diagrammes de Hasse sont respectivement, dans la figure~\ref{5posets}, $\poset_1$  et~$\poset_4$ sans l'élément~$w$.
Remarquons que l'équivalence est fausse pour les chaînes de Markov  à valeurs dans ces deux posets
(cf. les exemples~$1.1$
et~$4.5$ de~\cite{FillMachida}).

Bien plus, nous établissons la Proposition suivante
qui se fonde sur une technique permettant de \og relever\fg{} les générateurs monotones non complétement monotones
à valeurs dans un poset~$S'$ (par exemple de la figure~\ref{5posets}) en
un générateur monotone non complétement monotone
à valeurs dans un poset~$S$ contenant $S'$.
\medskip

{\bf Proposition~\ref{extension}}
Si un poset~$S$ admet comme sous-poset (cf. Définition~\ref{sub-poset})
 un poset~$S'$ dont le diagramme de Hasse est représenté dans les figures~\ref{5posets}
et~\ref{6posets}, alors les deux notions de monotonie ne sont pas équivalentes sur~$S$.
\medskip

Grâce à ce résultat et à des calculs directs, nous établissons également une liste exhaustive des posets de cardinal six
où l'équivalence n'est pas vérifiée. Nous renvoyons à la publication~\cite{DPLM} pour plus de détails.
La figure~\ref{6posets} donne les exemples caractéristiques de cardinal six
pour lesquels l'équivalence est fausse.

\selectlanguage{english}
\section{Introduction} \label{int}
The use of Markov chains in simulation has raised a number of
questions concerning qualitative and quantitative features of random processes,
in particular in connection with mixing properties. Among the features that are
useful in the analysis of effectiveness of Markov Chain Monte Carlo (MCMC) algorithms,
 a relevant role is played by monotonicity. Two notions of monotonicity have  been proposed,
for Markov chains taking values in a {\em partially ordered set} $S$ ({\em poset} from now on).
To avoid measurability issues, that are not relevant for our purposes, we shall always assume $S$
to be finite. Moreover, all Markov chains are implicitely assumed to be time-homogeneous.
\bd{mon1}
A Markov chain $(\eta_t)$, $t \in \R^+$ or $t \in \Z^+$, on the poset $S$, with transition
probabilities $P_t(x,y) := P(\eta_t = x|\eta_0=y)$, is said to be {\em monotone} if for each
pair $y,z \in S$ with $y \leq z$ there exists a Markov chain $(X_t(y,z))$ on $S \times S$ such that
\begin{enumerate}[(i)]
\item $X_0(y,z) = (y,z)$ a.s.
\item Each component $(X^i_t(y,z))$, $i=1,2$ is a Markov chain on $S$ with transition probabilities $P_t(x,y)$.
\item for all $t \geq 0$, $X^1_t(y,z) \leq X^2_t(y,z)$  a.s.
\end{enumerate}
\ed
There are various equivalent formulation of monotonicity. For instance, defining the
 {\em transition operator} $T_t f(y) := \sum_{x \in S} f(x) P_t(x,y)$, then the chain
is monotone if and only if $T_t$ maps increasing functions into increasing functions.
 This characterization can be turned (see Proposition~\ref{mon-cone}) into a simple algorithm
for checking monotonicity of Markov chains in terms of the element of the transition
matrix (in discrete-time) or in terms of the infinitesimal generator (in continuous-time).

For various purposes, including simulation, a stronger notion of monotonicity has been introduced.
\bd{mon2}
A Markov chain $(\eta_t)$, $t \in \R^+$ or $t \in \Z^+$, on the poset $S$, with transition
probabilities $P_t(x,y) := P(\eta_t = x|\eta_0=y)$, is said to be {\em completely monotone}
 if there exists a Markov chain $(\xi_t(\cdot))$ on $S^S$ such that
\begin{enumerate}[(i)]
\item $\xi_0(y) = y$ a.s.
\item For every fixed $z \in S$, the process $(\xi_t(z))$ is a Markov chain with transition probabilities $P_t(x,y)$.
\item If $y \leq z$, then for every $t \geq 0$ we have $\xi_t(y) \leq \xi_t(z)$ a.s.
\end{enumerate}
\ed
In other words, complete monotonicity means that we can simultaneously couple in
an order preserving way all processes leaving any possible initial state. This
property becomes relevant when one aims at sampling from the stationary measure
 of a Markov chain using the Propp \& Wilson algorithm (see \cite{ProppWilson}).

\noindent
If the transition probabilities, or the infinitesimal generator, are known, no simple
 general rule for checking complete monotonicity is known. Since, obviously, complete
 monotonicity implies monotonicity, a natural question is to determine for which posets
the converse is true. This problem has been completely solved in \cite{FillMachida} for
discrete-time Markov chain: in this case the following result holds.
\bt{fill-machida}
Every monotone Markov chain in the poset $S$ is also completely monotone if and only
if $S$ is {\em acyclic}, i.e. there is no loop $x_0,x_1, \ldots, x_n,x_{n+1} = x_0$ such that, for $i=0,1,\ldots,n$
\begin{enumerate}[(i)]
\item $\forall j=0,1,\ldots,n,\ j\neq i,\ x_i \neq x_{j}$,
\item Either $x_i < x_{i+1}$ or $x_i > x_{i+1}$,
\item $x_i \leq y \leq x_{i+1}$ or $x_i \geq y \geq x_{i+1}$ implies $y=x_i$ or $y=x_{i+1}$.
\end{enumerate}
\et
Our aim in this paper is to deal with the same problem in continuous-time, for {\em regular}
Markov chains, i.e. Markov chains possessing an infinitesimal
generator (or, equivalently, jumping a.s. finitely many times in any bounded time interval).
We have not been able to provide a complete link between discrete and continuous-time. It turns
out that if in a poset~$S$ monotonicity implies complete monotonicity in discrete-time, then the
same holds true in continous-time (see~Proposition~\ref{impl}). The converse is not true, however;
in the two four-points cyclic posets (the diamond and the bowtie, following the terminology in~\cite{FillMachida})
equivalence between monotonicity and complete monotonicity holds in continuous but not in discrete-time. There are,
however, five-points posets in which equivalence fails in continuous-time as well.

In this paper we do not achieve the goal of characterizing all posets for which equivalence
 holds. Via a computer-aided (but exact) method we give a complete list of five and six point
posets for which equivalence fails. Moreover we show that in each poset containing one of the
former as subposet, equivalence fails as well (this does not follow in a trivial way).

In Section \ref{pre} we give some preliminary notions, whose aim is to put the complete
monotonicity problem in continuous-time on a firm basis. Our main results are given in
Section \ref{equ}. All details not contained in this note will be given in the forthcoming paper~\cite{DPLM}.

\section{Preliminaries} \label{pre}

Let $(S,<)$ be a finite poset, and $L = (L_{x,y})_{x,y \in S}$ be the infinitesimal
generator of a regular Markov chain. Let $S_2 := S \times S \setminus \{(x,x): x \in S\}$.
 An infinitesimal generator is a matrix $L = (L_{x,y})_{x,y \in S}$ whose non-diagonal
elements are nonnegative, while the terms in the diagonal are given by $L_{x,x} = -\sum_{y \neq x} L_{x,y}$.
Thus $L$ may be identified with an element of the cone $(\R^+)^{S_2}$. Our main tool for understanding the
relations between the two notions of monotonicity consists in representing the sets of generators of
monotone and completely monotone Markov chains as subcones of $(\R^+)^{S_2}$.
\bd{upset}
A subset $\G \subseteq S$ is said to be an {\em up-set} if
$x \in \G \mbox{ and } x \leq y \ \Rightarrow \ y \in \G.$
\ed
The following proposition is a reformulation of a result in \cite{Massey87}.
By $\langle \cdot , \cdot \rangle$ we denote the Euclidean scalar product in $(\R)^{S_2}$.
\bp{mon-cone}
For  given  $x,y \in S$, $\G$ up-set, let $W^{\G,x,y} \in \R^{S_2}$ be defined by
\be{vectori_W}
W^{\G,x,y}_{v,z} = \left\{ \begin{array}{ll}
1 & \mbox{for } \begin{array}{l} x \leq y \not\in \G, \, v = y, \, z \in \G \mbox{ or} \\
						x \geq y \in \G, \, v=y, \, z \not\in \G
						\end{array} \\
-1 & \mbox{for } \begin{array}{l} x \leq y \not\in \G, \, v = x, \, z \in \G \mbox{ or} \\
						x \geq y \in \G, \, v=x, \, z \not\in \G
						\end{array} \\
0 & \mbox{in all other cases.}
\end{array}
\right.
\ee
Then the set ${\cal{G}}_{mon}$ of generators of monotone Markov chains is the subcone of the
elements $L \in (\R^+)^{S_2}$  satisfying the inequalities
\be{Cone_mon_H-representation}
\langle L, W^{\G,x,y} \rangle \geq 0 \ \ \mbox{for every } \G,x,y.
\ee
\ep
Let now ${\cal{M}}$ denote the set of increasing functions $S \ra S$.
For $f \in {\cal{M}}$, let $\mathbb{I}_f \in (\R^+)^{S_2}$ be defined by
\be{vectori_I}
(\mathbb{I}_f)_{x,y}  = \left\{ \begin{array}{ll} 1 & \mbox{if } f(x)=y \\ 0 & \mbox{otherwise.} \end{array}\right.
\ee
Our first result is the following description of the set of generators of completely
monotone Markov chains. Its proof is a rather straightforward consequence of the existence
of the {\em completely monotone coupling} (the process $(\xi_t)$ in Definition \ref{mon2}).
\bp{prop2.3}
The set ${\cal{G}}_{c.mon}$ of generators of completely monotone Markov chains is the cone
given by linear combination with nonnegative coefficients of the vectors $\mathbb{I}_f $, i.e.
\be{Cone_compl_mon_V-representation}
L = \sum_{f \in {\cal{M}}} \L_f \mathbb{I}_f, \textrm{ with } \L_f \geq 0.
\ee
\ep
Note that, for each $f \in {\cal{M}}$, $\G$ up-set, $x,y \in S$, we have
$ \langle \mathbb{I}_f,W^{\G,x,y} \rangle \geq 0,$
i.e. we recover the inclusion ${\cal{G}}_{c.mon}\subseteq {\cal{G}}_{mon}$. Our aim
is to determine for which posets the converse inclusion holds true.

A partial result in this direction is given by the following proposition, which
establishes a comparison with the discrete-time case. Its proof relies on analogous
representations in terms of cones for discrete-time transition matrix.
 \bp{impl}
 Suppose that in the poset $S$ monotonicity and complete monotonicity
are equivalent notions for discrete-time Markov chains. Then the equivalence holds
in continuous-time as well.
 \ep

\section{Extremal generators of monotone Markov chains: the monotonicity equivalence for "small" posets and extension to larger posets}
\label{equ}

According to the previous subsection, the equivalence between the complete monotonicity and the monotonicity of
 any Markov Chain on a poset~$S$ is equivalent to
\be{egualita_coni} {\cal{G}}_{c.mon}={\cal{G}}_{mon}. \ee
 In this section we answer the question for posets with a small cardinality.

 First note that the cases $\sharp S=2$,  $\sharp S=3$ are obvious:  in these cases $S$ is
linearly totally ordered (following the terminology of~\cite{FillMachida}). According to
Theorem~\ref{fill-machida}, there is equivalence for all the
discrete-time Markov chains and using the result of Proposition~\ref{impl}
the equivalence holds for continuous-time Markov chains as well.

In order to further investigate the equality~(\ref{egualita_coni}) we developed computer computations.
The cone~${\cal{G}}_{mon}$ is defined as intersection of half spaces in~(\ref{Cone_mon_H-representation}) (so called
\emph{H-representation}).
The cone~${\cal{G}}_{c.mon}$ is defined   by its extremal rays
in~(\ref{Cone_compl_mon_V-representation}) (so called
\emph{V-representation}). The software \emph{cdd+} (see \cite{cdd})
is able to compute exactly
one representation given the other one. This is a \emph{C++} implementation for convex polyhedron of
the Double Description Method (see for instance~\cite{FukudaProdon}).
Finding the extremal rays of the cone~${\cal{G}}_{mon}$ and the (minimal) set of inequalities defining the
 cone~${\cal{G}}_{c.mon}$, we are then able to know if ${\cal{G}}_{mon}\subseteq  {\cal{G}}_{c.mon}$ holds.

We operated by first using the software \emph{GAP} (see~\cite{GAP}) in order to
\begin{enumerate}[(i)]
\item  find the up-sets~$\G$ related to the poset~$S$, the
vectors~$W^{\G,x,y} \in \R^{S_2}$ and then identify the \mbox{H-representation} of~${\cal{G}}_{mon}$;
\item  \label{det_f_crescenti} %
compute all the increasing functions~$f \in {\cal{M}}$, identify the vectors~$\mathbb{I}_f \in (\R^+)^{S_2}$  and
then find the \mbox{V-representation} of~${\cal{G}}_{c.mon}$.
\end{enumerate}
We then use the software~\emph{cdd+} to produce the other representations of the cones, and the
software~\emph{Scilab} (see~\cite{Scilab})
to test if ${\cal{G}}_{mon}\subseteq  {\cal{G}}_{c.mon}$.

The difficulty in applying this method to posets with a high cardinality relies on the combinatorial complexity
of the step~(\ref{det_f_crescenti}) and on the time the software~\emph{cdd+} needs. Rather than to
$\sharp S$, this time is related to the number of facets
 of the cones, which comes from the partial order structure.
Nevertheless, we were able to completely study the cases when $\sharp S\leq 6$.

Here we briefly give the cases $\sharp S\leq 6$. For $\sharp S\geq 6$, the result of Proposition~\ref{extension}
gives the answer for some posets.
For a systematic analysis,  we refer to the
forthcoming work~\cite{DPLM}.

For $\sharp S=4$ the two relevant poset-structure are the diamond and the bowtie.
 Their Hasse-Diagram are, in Figure~\ref{5posets},
 $\poset_1$  and~$\poset_4$ respectively, without the node~$w$.
For those two posets, the computation procedure gives us that ${\cal{G}}_{mon}={\cal{G}}_{c.mon}$  holds.
Note that this result is known to be false in discrete-time,
 see for instance examples~$1.1$ and~$4.5$ in~\cite{FillMachida}.

For $\sharp S=5$, we can then state:
\bp{Prop3.1}
{The} only posets~$S$ with $\sharp S\leq 5$ such that (\ref{egualita_coni}) does not
hold are those whose Hasse-Diagram are presented in
Figure~\ref{5posets}.
\ep


\newsavebox{\Hasseprimo}
\begin{lrbox}{\Hasseprimo}{ %
 \xymatrix{ %
&  & \node \lulab{d} & &\\
& \node \lulab{b} \edge{ur} \edge{dr} &  & \node \rulab{c} \edge{ul} \edge{dl}&\\
&  & \node \ldlab{a}\edge{d} & \\
&  & \node \ldlab{w} & & \\
& & \ldlab{(\poset_1)} & & \\ } }
\end{lrbox}

\newsavebox{\Hassesecondo}
\sbox{\Hassesecondo}{\xymatrix{
 &  & \node \lulab{d}\\
 & \node \lulab{b} \edge{ur} \edge{dr}
 & \node \lulab{w} \edge{u}  \edge{d}
 & \node \rulab{c}\edge{ul} \edge{dl} & \\
 &  & \node \ldlab{a}\\
 & & &\\
 & & \ldlab{(\poset_2)}
&\\
}}

\newsavebox{\Hassetertio}
\sbox{\Hassetertio}{\xymatrix{
 &  & \node \lulab{d}\\
 & \node\lulab{b}\ar@{-}[ddrr]\edge{ur}& &\node \rulab{w} \edge{ul}\edge{dr}&\\
 &  &  & & \node \rulab{c} \edge{ul}  \ar@{-}[dl]&\\
 &  & &\node \ldlab{a}\\
 & & & \ldlab{(\poset_3)}& \\
} }

\newsavebox{\Hassequarto}
\sbox{\Hassequarto}{\xymatrix{
& \\
 & \node \lulab{c}\edge{d} \ar@{-}[ddrr] & &\node \rulab{d}\edge{dd}
\ar@{-}[ddll] &\\
 & \node \lulab{w}\edge{d}\\
 & \node \ldlab{a} & &\node\rdlab{b}&\\
 & & \ldlab{(\poset_4)}&\\
}}

\newsavebox{\Hassequinto}
\sbox{\Hassequinto}{ \xymatrix{
& & \node \lulab{w}\edge{dr}\edge{dl}\\
 & \node \lulab{c}\edge{dd} \ar@{-}[ddrr] & &\node \rulab{d}\edge{dd}
\ar@{-}[ddll] &\\
 & \\
 & \node \ldlab{a} & &\node\rdlab{b}&\\
 & & \ldlab{(\poset_5)}&
}}

\newsavebox{\HasseDoubleBowTie}
\sbox{\HasseDoubleBowTie}{
\xymatrix{
& \node \lulab{f} \edge{d} \edge{dr} & \node \lulab{e} & \node \rulab{d} \edge{dl} \edge{d}& \\
& \node \ldlab{a} \edge{ur}  & \node \ldlab{b} \edge{u} & \node \rdlab{c} \edge{ul} & \\
}}

\newsavebox{\HasseKrown}
\sbox{\HasseKrown}{\xymatrix{
& \node \lulab{f} \edge{d} \edge{dr} & \node \lulab{e} & \node \rulab{d} \edge{d}& \\
& \node \ldlab{a}  \ar@{-}[urr]  & \node \ldlab{b} \edge{u} & \node \rdlab{c} \edge{ul} & \\
}}

\newsavebox{\HasseKrownVariation}
\sbox{\HasseKrownVariation}{
\xymatrix{
& \node \lulab{f} \edge{d} \edge{dr} & \node \lulab{e} & \node \rulab{d} \edge{dl} \edge{d}& \\
& \node \ldlab{a} \edge{ur} \edge{urr} & \node \ldlab{b} \edge{u} & \node \rdlab{c} \edge{ul} & \\
}}

\newsavebox{\HasseKrownVariationBis}
\sbox{\HasseKrownVariationBis}{
\xymatrix{
& \node \lulab{f} \edge{d} \edge{dr} & \node \lulab{e} & \node \rulab{d} \edge{dl} \edge{d}& \\
& \node \ldlab{a} \edge{ur} \edge{urr} & \node \ldlab{b} \edge{u} & \node \rdlab{c} \edge{ul} \edge{ull} & \\
}}

\newsavebox{\HassePesce}
\sbox{\HassePesce}{\xymatrix{
& &   \node \lulab{d} & & \node \lulab{f} \edge{dl} & \\
& \node \lulab{b} \edge{ur} \edge{dr} &  & \node \ar@{}[rd]^<<{c} \edge{ul} \edge{dl}& & \\
& & \node \ldlab{a} & & \node \ldlab{e} \edge{lu} & \\
}}

\newsavebox{\HassePapillon}
\sbox{\HassePapillon}{\xymatrix{
& & \node \lulab{e} \edge{dl} \edge{dr} & & \node \lulab{f} \edge{dl} \edge{dr}& & \\
& \node \lulab{b} & & \node \node \ar@{}[rd]^<<{c} & & \node \rulab{d} & \\
& & & \node \rdlab{a} \ar@{-}[ull] \edge{u} \ar@{-}[urr] & & & \\
}}

\begin{figure}[h]
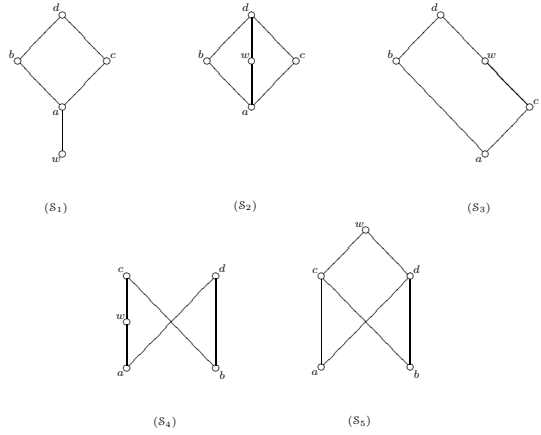

\centerline{ \scalebox{0.6}{\usebox{\Hasseprimo} %
\usebox{\Hassesecondo} %
\usebox{\Hassetertio}}}

\centerline{ \scalebox{0.6}{\usebox{\Hassequarto} %
\usebox{\Hassequinto} }}
\caption{\label{5posets}The five-points posets, for which there is no equivalence between the two notions of monotonicity.}
\end{figure}

\begin{figure}[h]
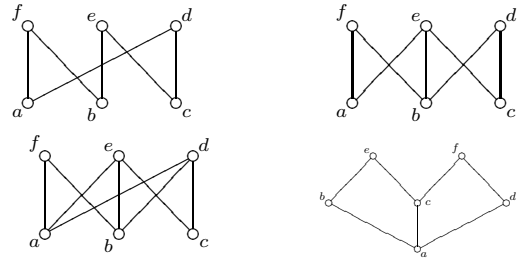

\centerline{ %
\usebox{\HasseKrown} %
\usebox{\HasseDoubleBowTie}}
\centerline{\usebox{\HasseKrownVariation}
 \scalebox{0.6}{\usebox{\HassePapillon}} }
\caption{\label{6posets}Some characteristic %
six-points posets, where there is no equivalence between the two notions of monotonicity.}
\end{figure}

We studied all five-points posets which are not linearly totally ordered.
For the posets previously mentioned, we found extremal rays $L=(L_{x,y})_{(x,y)\in S_2}$
of~${\cal{G}}_{c.mon}$ which are not in~${\cal{G}}_{mon}$.
One example for each poset is given in the figure~\ref{contro_esempi}. The non-mentioned components are equal to~$0$.

\begin{figure}[h]
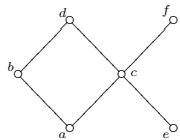

\centerline{\scalebox{0.7}{\usebox{\HassePesce}}}
\caption{\label{6posets-beispiel}Six-point poset, where there is no equivalence, since it admits $\poset_4$ as subposet.}
\end{figure}

Furthermore we briefly mention the following result, whose proof is to be found in~\cite{DPLM}.
\bd{sub-poset}
A subset~$S'$ of~$S$ is said to be an (induced) \impt{subposet} if for all $x,y\in S'$, $x\leq y$ in $S'$ is \impt{equivalent}
to $x\leq y$ in $S$.
\ed
\bp{extension}
If a poset~$S$ admits as subposet a poset~$S'$, whose Hasse-Diagram is presented in the figures~\ref{5posets}
and~\ref{6posets} then the monotonicity equivalence
({\it i.e. the equality~(\ref{egualita_coni})})  fails in~$S$ as well.
\ep
For instance, the poset of figure~\ref{6posets-beispiel} admits the poset~$\poset_4$ as subposet (consider it
without the node~$c$) and as consequence the equivalence of the two monotonicity does not hold on it.

\bp{six-posets} For $\sharp S=6$, the posets where the equivalence fails
are either the ones of figure~\ref{6posets} or the ones which admits as subposet a poset of figure~\ref{5posets}.
\ep

\begin{figure}[h]
$$\begin{array}{clcl}
\textrm{on }\poset_1, & L_1 & : & L_{a,w}= L_{b,w}=L_{c,w}=L_{d,b}=L_{d,c}=1 \\
\textrm{on }\poset_2, & L_2 & : & L_{a,b}=L_{d,b}=L_{w,b}=L_{c,d}=L_{c,a}=1 \\
\textrm{on }\poset_3, & L_3 & : & L_{d,c}=L_{b,c}=L_{a,c}=L_{w,a}=L_{c,a}=1 \\
\textrm{on }\poset_4, & L_4 & : & L_{a,b}=L_{w,b}=L_{c,d}=L_{d,b}=L_{b,d}=1 \\
\textrm{on }\poset_5, & L_5 & : & L_{w,c}=L_{w,d}=L_{c,a}=L_{d,a}=L_{b,a}=1 \\
\end{array}$$
\caption{\label{contro_esempi} Examples of monotone generators which are not complete monotone.}
\end{figure}

\section*{Acknowledgements}
P.-Y. Louis thanks the Mathematics Department of the University of Padova for financing stays at Padova.

\bibliographystyle{plain}

\end{document}